\documentclass[amscd,amssymb,11pt,reqno]{amsart}

\newtheorem{Thm}{Theorem}[section]

\newtheorem{Lem}[Thm]{Lemma}

\usepackage{graphicx}
\usepackage{latexsym}
\usepackage{amssymb}
\begin{document}


\vspace{1.5 cm}


\title[On projections of convex bodies]
      {Modified Shephard's problem \\ on projections of convex bodies}

\author{Vladyslav Yaskin}

\address{V.~Yaskin, Department of Mathematics, University of Oklahoma, Norman, OK 73019, USA.}
\email{vyaskin@math.ou.edu}


\begin{abstract}
We disprove a conjecture of A.~Koldobsky asking whether it is enough to compare $(n-2)$-derivatives
of the projection functions of two symmetric convex bodies in the Shephard problem in order to get
a positive answer in all dimensions.
\end{abstract}



\maketitle

\section{Introduction} Sections and projections of convex bodies have been actively studied for
many years. Although their properties exhibit certain duality, there is no clear picture explaining
this. The following two famous problems ask similar questions about sections and projection. Let
$K$ and $L$ be origin-symmetric convex bodies in $\mathbb{R}^n$. The Busemann-Petty problem asks
whether the assumption that all central hyperplane sections of $L$ have smaller volume than those
$L$ implies that $K$ has smaller $n$-dimensional volume. Its counterpart for projections is known
as the Shephard problem. It asks whether
$$\mathrm{vol}_{n-1}(K|\theta^\perp)\le \mathrm{vol}_{n-1}(L|\theta^\perp)$$ for all $\theta\in S^{n-1}$
implies that $$\mathrm{vol}_n(K)\le \mathrm{vol}_n(L).$$ The latter problem was solved
independently by Petty\cite{Pe} and Schneider\cite{S1}, who showed that the implication is correct
only if $n=2$. The solution to the Busemann-Petty problem was settled through the efforts of many
people (for historical details see \cite[pp.3-7]{K-book}) and it turned out that the answer to this
problem is affirmative only in dimensions $n\le 4$.

A unified approach to these problems was given by Koldobsky, Ryabogin and Zvavitch \cite{KRZ1},
\cite{KRZ2}, see also \cite[Section 5.1 and Section 8.4]{K-book}. They showed that these two
problems are essentially of the same nature, if treated with the help of Fourier analysis.

Koldobsky \cite{K-AdvAppl} and  Koldobsky, Yaskin, Yaskina \cite{KYY} considered a modification of
the Busemann-Petty problem, which gave a positive answer to the problem in all dimensions. Namely,
for an origin-symmetric convex body $K$ in $\mathbb{R}^n$, define the section function
$$
S_K(\xi)=\mbox{vol}_{n-1}(K\cap\xi^\perp), \qquad \xi\in S^{n-1},
$$
where $\xi^\perp $ is the central hyperplane in $\mathbb{R}^n$ orthogonal to $\xi$, and extend
$S_K$ from the sphere to the whole $\mathbb{R}^n$ as a homogeneous function of degree $-1$. Let
$\Delta$ be the Laplace operator on $\mathbb{R}^n$. It was proved that for two origin-symmetric
infinitely smooth convex bodies $K,L$ in $\mathbb R^n$   and $\alpha\in \mathbb R,\ \alpha \ge
n-4,$ the condition
\begin{eqnarray}\label{eqn:condition1}
(-\Delta)^{\alpha/2} S_K(\xi)\le (-\Delta)^{\alpha/2} S_L(\xi), \qquad \forall \xi\in S^{n-1}
\end{eqnarray}
implies that $\mathrm{vol}_n(K)\le \mathrm{vol}_n(L),$ while for $\alpha<n-4$ this is not
necessarily true.

Koldobsky conjectured that a similar result must hold for projections, with $\alpha = n-2$ being
the critical value after which the answer becomes affirmative. This conjecture was based on the
following. First of all, for all origin-symmetric convex bodies $\|x\|_K^{-n+3}$ is a positive
definite distribution, see \cite[Section 4.2]{K-book}. In addition the authors of \cite{KYY} showed
that $\|x\|_K^{-1} |x|_2^{-n+4}$ is also positive definite, which corresponds to the borderline
case of the result mentioned above. (In general $\|x\|_K^{-n+p+3}\cdot |x|^{-p}_2$ is positive
definite for a certain range of $p$). One can see that all these functions have a common property:
they are homogeneous of degree $-n+3$, and therefore it seemed plausible that $\|x\|_K
|x|_2^{-n+2}$ should also be positive definite. (This would correspond to the case $p=n-2$ above).
However in \cite{KYY} the authors were unable to extend the proof to this case.

Here we prove that $\|x\|_K |x|_2^{-n+2}$ is not necessarily positive definite, giving a negative
answer to the conjecture of Koldobsky. This seems to be one of not many examples where the direct
analogy between sections and projections does not hold.

For other generalizations of the Shephard problem see \cite{B}, \cite{GZ}, \cite{RZ}.

{\bf Acknowledgments.} The author would like to thank Professors Paul Goodey and Alexander
Koldobsky for fruitful discussions.

\section{Convex Geometry and the Fourier transform}

The standard references here are the books by Gardner \cite{G}, Schneider \cite{S} and Koldobsky
\cite{K-book}. Let $K$ be an origin-symmetric star body in $\mathbb{R}^n$. The {\it Minkowski
functional} of $K$ is defined as
$$
\|x\|_K=\min \{a\ge 0: x \in a K \}.
$$
The function $\rho_K(x)=\|x\|_K^{-1}$ is called  the {\it radial function} of $K$. If $x\in
S^{n-1}$, $\rho_K(x)$ is the distance from the origin to the boundary of $K$ in the direction of
$x$.


We say that a body $K$ is infinitely smooth if its radial function $\rho_K$ restricted to the unit
sphere $S^{n-1}$ belongs to the space $C^\infty(S^{n-1})$ of infinitely differentiable functions on
the unit sphere.

Throughout the paper we use the Fourier transform of distributions. The Fourier transform of a
distribution $f$ is defined by $\langle\hat{f}, \phi\rangle= \langle f, \hat{\phi} \rangle$ for
every test function $\phi$ from the Schwartz space $ \mathcal{S}$ of rapidly decreasing infinitely
differentiable functions on $\mathbb R^n$. For any even distribution $f$, we have $(\hat{f})^\wedge
= (2\pi)^n f$.

In particular we are interested in  the Fourier transform of homogeneous functions on
$\mathbb{R}^n$. We will need the following version of Parseval's formula on the sphere proved by
Koldobsky, see e.g. \cite[p.66]{K-book}.
\begin{Lem}\label{Lem:Parseval} If $K$ and $L$ are origin-symmetric
infinitely smooth star bodies in $\mathbb{R}^n$ and $0<p<n$, then $(\|x\|_K^{-p})^\wedge$ and
$(\|x\|_L^{-n+p})^\wedge$ are continuous functions on $S^{n-1}$ and
$$\int_{S^{n-1}} \left(\|x\|_K^{-p}\right)^\wedge (\xi) \left(\|x\|_L^{-n+p}\right)^\wedge
(\xi)d\xi= (2\pi)^n \int_{S^{n-1}} \|x\|_K^{-p}\|x\|_L^{-n+p} dx.$$
\end{Lem}

\noindent {\bf Remark.} The preceding   lemma  was formulated for Minkowski functionals, but in
fact it holds true for arbitrary infinitely differentiable even functions on the sphere extended to
$\mathbb{R}^n\setminus\{0\}$ as homogeneous functions of corresponding degrees.

We say that a distribution $f$ is {\it positive definite} if its Fourier transform is a positive
distribution, in the sense that $\langle \hat{f},\phi \rangle \ge 0$ for every non-negative test
function $\phi$.

The next result from \cite{GKS} will be our main tool in determining  whether a homogeneous
function represents a positive definite distribution, see also \cite[p.60]{K-book}.

\begin{Thm}\label{Thm:GKS}{\rm (Gardner, Koldobsky, Schlumprecht)}
Let $K$ be an infinitely smooth origin-symmetric star body in $\mathbb{R}^n$, and let
$k\in\mathbb{N}\cup  \{0\}$, $k\ne n-1$. Suppose that $\xi\in S^{n-1}$, and let $A_{K,\xi}$ be the
corresponding parallel section function of $K$: $A_{K,\xi}(z)=\int_{K\cap ( x,\xi)=z} dx$.
\newline

(a) If $q$ is not an integer, $-k-1<q<k$, then
\begin{eqnarray*}(\|x\|_K^{-n+q+1})^\wedge(\xi)=\frac{\pi(n-q-1)}{\Gamma(-q) \cos\frac{\pi
q}{2} }\times\hspace{5.5cm}\\
\times\int_0^\infty
\frac{A_{K,\xi}(z)-A_\xi(0)-{A}_{K,\xi}''(0)\frac{z^2}{2}-\cdots-A_{K,\xi}^{(k-1)}(0)\frac{z^{k-1}}{(k-1)!}}{z^{q+1}}dz.
\end{eqnarray*}

(b) If $k$ is an even integer, then
$$
(\|x\|_K^{-n+k+1})^\wedge(\xi)=(-1)^{k/2}\pi (n-k-1) A_{K,\xi}^{(k)}(0).
$$

(c) If $k$ is an odd integer, then
\begin{eqnarray*}
(\|x\|_K^{-n+k+1})^\wedge(\xi) =(-1)^{(k+1)/2}2(n-1-k)k!\times\hspace{3.5cm}\\
\times\int_0^\infty
\frac{A_{K,\xi}(z)-A_{K,\xi}(0)-{A}_{K,\xi}''(0)\frac{z^2}{2}-\cdots-A_{K,\xi}^{(k-1)}(0)\frac{z^{k-1}}{(k-1)!}}{z^{k+1}}dz.
\end{eqnarray*}
Here $A_{K,\xi}^{(k)}$ stands for the derivative of the order $k$ and the Fourier transform is
considered in the sense of distributions.
\end{Thm}
\noindent {\bf Remarks.} (i) The previous theorem implies that for infinitely smooth bodies the
Fourier transform of $\|x\|^{-n+q+1}$ restricted to the unit sphere  is a continuous function (see
also \cite[Section 3.3]{K-book}).

(ii) If $k=0$, then part (a) of the theorem reads as follows. For $-1<q<0$,
$$(\|x\|_K^{-n+q+1})^\wedge(\xi)=\frac{\pi(n-q-1)}{\Gamma(-q) \cos\frac{\pi q}{2} } \int_{S^{n-1}}
|(\theta,\xi)|^{-q-1} \|\theta\|_K^{-n+q+1} d\theta.$$ In particular, if $-1<q<0$, then
$(\|x\|_K^{-n+q+1})^\wedge $ is a non-negative function on the sphere for any star body $K$.

An extension of Theorem \ref{Thm:GKS} to the case when $k=n-1$ was given in \cite{KKYY}.
\begin{Thm}\label{Thm:n-1-deriv}
Let $K$ be an infinitely smooth origin-symmetric star body in $\mathbb{R}^n$. Extend
$A_{K,\xi}^{(n-1)}(0)$ to a homogeneous function of degree $-n$ of the variable $\xi \in
\mathbb{R}^n \setminus \{0\}$. Then $(\ln \|\cdot\|_K)^\wedge$ is a continuous function on
$\mathbb{R}^n\setminus \{0\}$ and
\begin{eqnarray}\label{eqn:A^{(n-1)}}
A_{K,\xi}^{(n-1)}(0)=-\frac{\cos(\pi (n-1)/2)}{\pi} \left( \ln \|\cdot\|_K\right)^\wedge (\xi),
\end{eqnarray}
as distributions (of the variable $\xi$) acting on test functions with compact support outside of
the origin. In particular,

i) if $n$ is odd

$$\left( \ln \|x\|_K\right)^\wedge (\xi)=(-1)^{(n+1)/2}\pi A_{K,\xi}^{(n-1)}(0),  \ \ \xi \in \mathbb{R}^n\setminus \{0\},$$

ii) if $n$ is even, then for    $\xi \in \mathbb{R}^n\setminus \{0\}$,

$$\left( \ln \|x\|_K\right)^\wedge (\xi)=\hspace{9cm}$$ $$=a_n\int_0^{\infty}\frac{A_{K,\xi}(z)-A_{K,\xi}(0)-
A_{K,\xi}''(0)\frac{z^2}{2}-...-A^{n-2}_{K,\xi}(z)\frac{z^{n-2}}{(n-2)!}}{z^n}dz,$$ where
$a_n=2(-1)^{n/2+1}(n-1)!$

\end{Thm}

The {\it support function} of a convex body $K$ in $\mathbb{R}^n$ is defined by
$$h_K(x)=\max_{\xi \in K}(x,\xi),\ \ x\in \mathbb{R}^n. $$
If $K$ is origin-symmetric, then $h_K$ is the Minkowski norm of the polar body $K^*$.

Let $\mathrm{vol}_{n-1}(K|\theta^\perp)$ denote the $(n-1)$-dimensional volume of the orthogonal
projection of the body $K$ onto the hyperplane orthogonal to $\theta$. The following is the
well-known Cauchy formula \cite[p. 361]{G}:
$$\mathrm{vol}_{n-1}(K|\theta^\perp)=\frac12 \int_{S^{n-1}} |(\xi,\theta)| dS_{n-1} (K,\xi),$$
where $dS_{n-1} (K,\xi)$ is the surface area measure of $K$ (\cite[p. 351]{G}). A convex body $K$
is said to have a {\it curvature function} $f_K$, if its surface area measure $dS_{n-1} (K,\xi)$ is
absolutely continuous with respect to Lebesgue measure $d\sigma_{n-1} $ on $S^{n-1}$ and
$$\frac{dS_{n-1}(K,\cdot)}{d\sigma} = f_K(\cdot) \in L^1(S^{n-1}).$$

If $K$ is an infinitely smooth body with positive curvature, then $f_K(\theta)$ is the reciprocal
of the Gauss curvature at the boundary point with unit normal $\theta$, see \cite[p. 419]{S}.
Abusing notations, we will also denote by $f_K$ the extension of $f_K$ to $\mathbb{R}^n$ as a
homogeneous function of degree $-n-1$. Koldobsky, Ryabogin and Zvavitch \cite{KRZ1} proved that if
a body $K$ has a curvature function, then \begin{equation}\label{Formula for proj}
\mathrm{vol}_{n-1}(K|\theta^\perp)= -\frac{1}{\pi} \widehat{f_K}(\theta) , \qquad \forall \theta
\in S^{n-1}.
\end{equation}

Let $\Delta$ be the Laplace operator on $\mathbb{R}^n$. The fractional powers of the Laplacian of a
distribution $g$ are defined by
\begin{equation}\label{Def:Laplacian}
 (-\Delta)^{\alpha/2}g =
\frac{1}{(2\pi)^n}( |x|_2^\alpha \hat{g}(x))^\wedge,
\end{equation}
where  the Fourier transform is considered in the sense of distributions, and $|x|_2$ stands for
the Euclidean norm in $\mathbb R^n.$ Using the connection between the Fourier transform and
differentiation, one can see that for an even integer $\alpha$ and an even distribution $g$ this
definition gives the standard Laplacian applied $\alpha/2$ times.


If $K$ is an infinitely smooth body with positive Gauss curvature, then $f_K$ is an infinitely
differentiable function on the sphere (because the Gauss curvature is the determinant of the
Weingarten map, which is infinitely differentiable and non-singular in our case, see
\cite[pp.104-109]{S}). Consider the projection function $\mathrm{vol}_{n-1}(K|(\cdot)^\perp)$ and
extend it from the sphere to $\mathbb{R}^n$ as a homogeneous function of degree 1. Using
(\ref{Formula for proj}) and (\ref{Def:Laplacian}), we get
\begin{equation}\label{eqn:Def-Laplace}
(-\Delta)^{\alpha/2}\mathrm{vol}_{n-1}(K|\theta^\perp) = -\frac{1}{\pi}( |x|_2^\alpha
f_K(x))^\wedge(\theta).
\end{equation}

Since $|\cdot|_2^\alpha f_K$ is  infinitely differentiable, $(-\Delta)^{\alpha/2}
\mathrm{vol}_{n-1}(K|(\cdot)^\perp)$ is a continuous function on the sphere (\cite[Lemma
3.16]{K-book} or remark (i) after Theorem \ref{Thm:GKS}).


\section{Main results}
Let us start with a result in the positive direction.
\begin{Thm}\label{Thm:main2}
Let $n \le \alpha < n+1$. Let  $K$, $L\subset \mathbb{R}^n$, $n\ge 3$, be infinitely smooth
origin-symmetric convex bodies with
positive Gauss curvature   such that 
 $$(-\Delta)^{\alpha/2}\mathrm{vol}_{n-1}(K|\theta^\perp)\ge (-\Delta)^{\alpha/2} \mathrm{vol}_{n-1}(L|\theta^\perp), \qquad \forall \theta\in S^{n-1}.$$
Then $\mathrm{vol}_n(K)\le \mathrm{vol}_n(L)$.

\end{Thm}

\noindent{\bf Proof.} Recall the following formula \cite[p. 354]{G}
$$\mbox{vol}_n(L)=\frac1n\int_{S^{n-1}}h_L(\theta)f_L(\theta)d\theta,$$
where $h_L$  and $f_L$ are the support function  and curvature function of the body $L$
correspondingly. Therefore using the fact that
$h_L=\|\cdot\|_{L^*}$ 
we get
\begin{eqnarray*}
\mbox{vol}_n(L)&=&\frac1n\int_{S^{n-1}}\|\theta\|_{L^*}f_L(\theta)d\theta\\
&=&\frac{1}{n}\int_{S^{n-1}}|\theta|_2^{-\alpha}\|\theta\|_{L^*}|\theta|_2^{\alpha}f_L(\theta)d\theta.
\end{eqnarray*}

Since $n\le \alpha < n+1$,   remark (ii) after Theorem \ref{Thm:GKS}  implies that the Fourier
transform of $|x|_2^{-\alpha}\|x\|_{L^*}$ is a non-negative function on the sphere. Applying the
spherical version of Parseval's formula (Lemma \ref{Lem:Parseval}), we get
\begin{eqnarray*}
&=&\frac{1}{(2\pi)^n n}\int_{S^{n-1}}(|x|_2^{-\alpha}\|x\|_{L^*})^\wedge(\xi)(|x|_2^{\alpha}f_L(x))^\wedge(\xi) d\xi \\
&=&-\frac{\pi}{(2\pi)^n n}\int_{S^{n-1}}(|x|_2^{-\alpha}\|x\|_{L^*})^\wedge(\xi) (-\Delta)^{\alpha/2}\mbox{vol}_{n-1}(L|\xi^\perp)d\xi \\
&\ge&-\frac{\pi}{(2\pi)^n n}\int_{S^{n-1}}(|x|_2^{-\alpha}\|x\|_{L^*})^\wedge(\xi) (-\Delta)^{\alpha/2}\mbox{vol}_{n-1}(K|\xi^\perp)d\xi \\
&=&\frac1n\int_{S^{n-1}}\|\theta\|_{L^*}f_K(\theta)d\theta=\frac1n\int_{S^{n-1}}h_L(\theta)f_K(\theta)d\theta\\
&=&V_1(K,L),
\end{eqnarray*}
where $V_1(K,L)$ is the mixed volume, also denoted by $V(K,...,K,L)$, see \cite[p.353]{G},
\cite[p.275]{S}.

Therefore we have $V_1(K,L)\le \mbox{vol}_n(L)$. Applying Minkowski's first inequality \cite[p.
317]{S} we get

$$\mbox{vol}_n(L)^\frac{1}{n} \mbox{vol}_n(K)^\frac{n-1}{n}\le V_1(K,L)\le \mbox{vol}_n(L),$$
and hence $$\mbox{vol}_n(K)\le \mbox{vol}_n(L).$$

\qed

  {\bf Remark.} Comparing   the previous theorem with the original Shephard problem, one
can observe that the inequality for the projections gets reversed. This happens because the answer
to Shephard's problem is affirmative if $L$ is a polar projection body, that is the Fourier
transform of $\|\cdot\|_{L^*}$ is a {\it negative} distribution outside of the origin, see
\cite[pp.155-160]{K-book}. On the other hand, as we have seen, if this norm is multiplied by  the
Euclidean norm to the appropriate power, then the Fourier transform of $|x|_2^{-\alpha}\|x\|_{L^*}$
becomes a {\it positive} distribution.

\begin{Lem}\label{Lem:main}

 Let $n-2 \le \alpha < n$, $\alpha\ne 1$. Then there exists an origin-symmetric convex body $L$ in
$\mathbb{R}^n$, $n\ge 3$, such that $|x|_2^{-\alpha}\|x\|_{L}$ is not a positive definite
distribution.

\end{Lem}
\noindent{\bf Proof.} 
First consider the case $n-2 < \alpha < n$. For a large $N>0$ let $L$ be an ellipsoid with the
norm:
$$\|x\|_L = (x_1^2+\cdots+x_{n-1}^2+ N x_n^2)^{1/2}.$$


Define a star body $K \subset \mathbb{R}^n$ by the formula:
$$\rho_K(\theta) = \rho_L^{\frac{1}{1 -\alpha }} (\theta), \qquad \theta \in S^{n-1},$$
where $\rho_K$ and $\rho_L$ are the radial functions of the bodies $K$ and $L$ correspondingly. One
can see that
$$|x|_2^{-\alpha}\|x\|_{L} =\left( |x|_2^{-\alpha/(1 - \alpha)}\|x\|_{L}^{1/( 1-\alpha )
}\right)^{-\alpha +1 } = \|x\|_K^{-\alpha+1}, \qquad \forall x \in \mathbb{R}^n\setminus\{0\}. $$

Using Theorem \ref{Thm:GKS} with $q=n-\alpha \in (0,2)$ we get
$$(\|x\|_K^{-\alpha+1})^\wedge(\xi)=\frac{\pi(\alpha-1)}{\Gamma(\alpha - n) \cos\frac{\pi
(n-\alpha)}{2} }\int_0^\infty t^{-n+\alpha-1}(A_{K,\xi}(t)-A_{K,\xi}(0))dt,$$ where the case
$n-\alpha=1$ is understood in the sense of part (c) of the aforementioned theorem.

Note that ${\Gamma(\alpha - n) \cos\frac{\pi (n-\alpha)}{2} } \le 0$ for $\alpha\in (n-2,n)$, where
$\alpha=n-1$ is again understood in terms of the limit, so we need to prove that for some $\xi$
\begin{equation}\label{int}
\int_0^\infty t^{-n+\alpha-1}(A_{K,\xi}(t)-A_{K,\xi}(0))dt > 0.
\end{equation}

Let $\xi$ be the direction of the $x_n$-axis. Let $[-t_0,t_0]$ be the support of $A_{K,\xi}(t)$,
then
$$\int_0^\infty t^{-n+\alpha-1}(A_{K,\xi}(t)-A_{K,\xi}(0))dt= $$
$$=\int_0^{t_0} t^{-n+\alpha-1}(A_{K,\xi}(t)-A_{K,\xi}(0))dt- \int_{t_0}^\infty t^{-n+\alpha-1}
A_{K,\xi}(0) dt$$
\begin{equation}\label{int2}=\int_0^{t_0}
t^{-n+\alpha-1}(A_{K,\xi}(t)-A_{K,\xi}(0))dt - \frac{A_{K,\xi}(0)}{n-\alpha }t_0^{-n+\alpha}.
\end{equation}



Introduce the following coordinates on the sphere $S^{n-1}$.  Every $\theta \in S^{n-1}$ can be
written as
$$\theta = \cos\phi \cdot \zeta  + \sin\phi \cdot \xi,$$
where $-\pi/2\le \phi \le \pi/2$ and $\zeta \in S^{n-1}\cap \xi^\perp$.

Since we are interested in the sections of $K$ perpendicular to $\xi$, its axis of revolution, by
abuse of notation we will denote by $\rho_K(\phi)$ the radial function of those $\theta \in
S^{n-1}$ that make an angle $\phi$ with the plane $\xi^\perp$. Explicitly it equals $$\rho_K(\phi)
= (\cos^2 \phi + N \sin^2 \phi )^{1/(2\alpha-2)}.$$

One can check that $t=\sin \phi \cdot \rho_K(\phi)$ is an increasing function of the angle $\phi\in
(0,\pi/2)$, therefore all the sections of $K$ by hyperplanes orthogonal to $\xi$ are
$(n-1)$-dimensional disks. Moreover, one can see that $t_0=N^{\frac{1}{2 \alpha - 2}}$, which
implies that the last term in (\ref{int2}) approaches zero as $N$ tends to infinity.

It will be more convenient to work with $\phi$ instead of $t$. $A_{K,\xi}$ as a function of $\phi$
looks as follows.
\begin{eqnarray*}
A_{K,\xi}(t(\phi)) &= &\omega_{n-1} (\cos \phi \cdot \rho_K(\phi))^{n-1}\\
&=& \omega_{n-1} (\cos \phi)^{n-1} ( \cos^2\phi +N \sin^2 \phi)^{\frac{n-1}{2\alpha -2}},
\end{eqnarray*} where $\omega_{n-1}$ is the volume of the unit $(n-1)$-dimensional Euclidean
ball.


Now consider the integral term from (\ref{int2}). Making change of the variable $t=\sin \phi \cdot
(\cos^2 \phi + N \sin^2 \phi )^{1/(2\alpha-2)}$, we get
$$\int_0^{t_0} t^{-n+\alpha-1}(A_{K,\xi}(t)-A_{K,\xi}(0))dt=$$
$$=\omega_{n-1} \int_0^{\pi/2} (\sin \phi)^{-n+\alpha-1} ( \cos^2\phi +N \sin^2
\phi)^{\frac{-n+\alpha-1}{2\alpha -2 }} \times$$ $$\times \left(  (\cos \phi)^{n-1} ( \cos^2\phi +N
\sin^2 \phi)^{\frac{n-1}{2\alpha -2 }} - 1\right)\times$$
\begin{equation}\label{Eqn:Integral}
\times ( \cos^2\phi +N \sin^2 \phi)^{\frac{1}{2\alpha -2 }-1}(\cos^3 \phi
+(N+\frac{N-1}{\alpha-1})\cos\phi \sin^2\phi) d\phi.
\end{equation}

Now we want to find the intervals where the integrand is positive or negative. So we need to solve
the equation
$$(\cos \phi)^{n-1} ( \cos^2\phi +N
\sin^2 \phi)^{\frac{n-1}{2\alpha -2 }} - 1=0,$$ which is equivalent to
\begin{equation}\label{Eqn:phi}
(\cos \phi)^{2\alpha} + N (\cos \phi)^{2\alpha - 2} \sin^2 \phi =1.
\end{equation}

By showing that the function in the left hand side is first increasing and then decreasing to zero,
one can see that the equation has two roots on the interval $[0, \pi/2]$. One root is obvious:
$\phi_1 =0$. In order to determine the second root $\phi_2$, note that the maximum of the function
in question is achieved when $\phi$ is roughly $\arccos \sqrt{1 - 1/\alpha}$, assuming $N$ is
large. Therefore   (\ref{Eqn:phi}) together with the inequality $\phi_2   \gtrsim \arccos \sqrt{1 -
1/\alpha}$ gives
$$N (\cos \phi_2)^{2\alpha - 2} \le C(\alpha),$$ and hence $\phi_2=\pi/2 -
o(N^{\frac{1}{2-2\alpha}})$.

Now break the integral (\ref{Eqn:Integral}) into two parts according to where the integrand is
positive or negative. It is negative on the interval $(\pi/2 - o(N^{\frac{1}{2-2\alpha}}), \pi/2)$
and one can easily show that the absolute value of the integral here is bounded above by
$$C N^{\frac{\alpha - n -2}{2\alpha -2}},$$
which approaches zero as $N$ tends to infinity.

In order to estimate from below the positive part of the integral (\ref{Eqn:Integral}) it is enough
to consider the interval $[\pi/4, \pi/3]$. One can check that when $N$ is large, the integral has
order
$$C N^{\frac{1}{2}},$$ which approaches infinity as $N$ gets large.
The inequality (\ref{int}) follows.

Now consider $\alpha = n-2$. In this case Theorem \ref{Thm:GKS} gives $$
(\|x\|_K^{-\alpha+1})^\wedge(\xi)= {\pi(1 -\alpha)}A_{K,\xi}^{''}(0)<0.$$ The latter inequality
follows by direct computation.

\qed

The previous Lemma says nothing about the case when $\alpha=1$ (and therefore $n=3$). It may seem
that the right analog would be to analyze the sign of  $(|x|_2^{-1} \|x\|_L  + \ln|x|_2)^\wedge$.
But in fact, as one will see later,  the following result is needed.

\begin{Lem}\label{Lem:alpha=1}
 There exists an origin-symmetric convex body $L$ in $\mathbb{R}^3$, such that the Fourier
transform of $$  |x|_2^{-1}\|x\|_{L} - \frac{\int_{S^2} \|\theta\|_{L}d\theta
}{4\pi(1+\Gamma'(1))}\ \ln |x|_2
$$ is not a positive distribution outside of the origin. Here $\Gamma'$ is the derivative of the
Gamma-function.
\end{Lem}

\noindent{\bf Proof.}
For $N>0$ large enough consider the following planar curve defined in polar coordinates by
 \begin{eqnarray}\label{eqn:curve}
 \rho(\phi)=\cos^N \phi.
 \end{eqnarray}
 Take only that part of the curve where the angle $\phi$ belongs to the interval
 $$[-\arcsin\frac{1}{\sqrt{N+1}}, \arcsin\frac{1}{\sqrt{N+1}}].$$ One can check that the
 end-points of the interval correspond to the extreme values of the $y$-coordinate (altitude). Rotate this arc around the $y$-axis and then
 attach two disks at the top and bottom to get a closed surface. Denote by $L$  the convex body bounded by this surface. Consider a star body $K$ given by the formula
$$ \|\theta\|_K  = \exp({\|\theta\|_L }), \qquad \theta \in S^{2}.$$
Therefore the radial function of $K$ equals $$\rho_K(\theta) = \exp(-\rho_L^{-1}(\theta)), \qquad
\theta \in S^{2}.$$
One can also see that 
$$\ln \|x\|_K = |x|_2^{-1} \|x\|_L  + \ln|x|_2, \qquad x \in \mathbb{R}^3\setminus\{0\}.$$
Let $\xi$ be the direction of the axis of revolution of $L$. Since $n=3$, by Theorem
\ref{Thm:n-1-deriv} we have
$$\left( \ln \|x\|_K\right)^\wedge (\xi)= \pi A_{K,\xi}^{''}(0).$$
(In fact this formula can only be applied if the body is smooth enough, but let us ignore this
problem for a while and address it at the end of the proof).

If  we denote by $\rho_K(\phi)$ the radial function of those $\theta \in S^{2}$ that make an angle
$\phi$ with the plane $\xi^\perp$, then
 $$A_{K,\xi}(\phi) = \pi (\cos\phi \cdot \rho_K(\phi))^2 = \pi \cos^2\phi\cdot \exp(-2\rho_L^{-1}(\phi)).$$
 Using that for small $\phi$ the function $\rho_L$ is given by formula (\ref{eqn:curve}) we get
 \begin{equation}\label{A''}
 A_{K,\xi}^{''}(0)=-2\pi e^{-2}(N-1).
 \end{equation}

On the other hand from the construction of the body $L$ it follows that  $L$ has smallest radius in
the direction of $\xi$. Therefore for all $\theta \in S^2$,  $$\rho_L(\theta)\ge \rho_L(\xi) =
\left(\cos^N\phi \sin\phi\right)|_{\phi=\arcsin\frac{1}{\sqrt{N+1}}} \simeq
\frac{C_1}{\sqrt{N+1}},$$ which implies
$$\int_{S^2} \|\theta\|_{L}d\theta \le C_2 \sqrt{N+1},$$ for some constants $C_1$ , $C_2>0$.

Also notice that part (i) of Theorem \ref{Thm:n-1-deriv} gives  $$(\ln |x|_2)^\wedge(\theta) =
-2\pi^2, \qquad \forall \theta \in S^2.$$  Therefore we have

$$  \left(|x|_2^{-1}\|x\|_{L} - \frac{\int_{S^2} \|\theta\|_{L}d\theta }{4\pi(1+\Gamma'(1))}\ \ln
|x|_2\right)^\wedge (\xi)
$$
$$=\left( \ln \|x\|_K - \frac{1+\Gamma'(1)+\frac{1}{4\pi}\int_{S^2} \|\theta\|_{L}d\theta }{1+\Gamma'(1)}\ \ln
|x|_2\right)^\wedge (\xi) $$
$$
=-2\pi e^{-2}(N-1)+2\pi^2\frac{1+\Gamma'(1)+\frac{1}{4\pi}\int_{S^2} \|\theta\|_{L}d\theta
}{1+\Gamma'(1)}$$
$$
\le -2\pi e^{-2}(N-1)+C \sqrt{N+1} <0,$$ for $N>0$ large enough.

Formally the above computations are not quite legitimate since $L$ is not infinitely smooth. But
one can approximate $L$ by an origin-symmetric infinitely smooth convex body without loosing the
sign in the last inequality. Specifically, one has to smooth out the body in a small neighborhood
of ${\phi=\arcsin\frac{1}{\sqrt{N+1}}}$. This operation will not affect (\ref{A''}). On the other
hand one can also assure that $\int_{S^2} \|\theta\|_{L}d\theta$ does not change much.

  \qed

\begin{Thm}\label{Thm:main3}
Let $n-2 \le \alpha < n$. There are convex origin-symmetric bodies $K,L\subset \mathbb{R}^n$, $n\ge
3$ such that
\begin{eqnarray}\label{eqn:Laplacians}
(-\Delta)^{\alpha/2}\mathrm{vol}_{n-1}(L|\theta^\perp)\le
(-\Delta)^{\alpha/2}\mathrm{vol}_{n-1}(K|\theta^\perp), \qquad \forall \theta\in S^{n-1},
\end{eqnarray}
but $$\mathrm{vol}_n(L)< \mathrm{vol}_n(K).$$
\end{Thm}

\noindent {\bf Proof.} First assume that $\alpha\ne 1$. Lemma \ref{Lem:main} guarantees that there
exists an ellipsoid $K^*$, such that $(|x|_2^{-\alpha} \|x\|_{K^*})^\wedge(\xi)< 0$ for some
direction $\xi$. Let $K$ be the polar body of $K^*$. Since $K$ is again an ellipsoid,   its
curvature function $f_K$ is well-defined.

Let $\Omega=\{\theta \in S^{n-1}:(|x|_2^{-\alpha} \|x\|_{K^*})^\wedge(\theta)< 0\}$ and let $v\in
C^\infty (S^{n-1})$ be a non-negative even function supported in $\Omega$. Extend $v$ to a
homogeneous function $|x|_2^{1-\alpha}v(x/|x|_2)$ of degree $1- \alpha$ on $\mathbb{R}^n$. By
\cite[Lemma 5]{K-Israel99} the Fourier transform of $|x|_2^{1-\alpha}v(x/|x|_2)$ is equal to
$|x|_2^{-n-1+\alpha}g(x/|x|_2)$ for some function $g\in C^\infty(S^{n-1})$. Choose an
$\varepsilon>0$ small enough and define
$$ f_L(x)=f_K(x) + \varepsilon |x|_2^{-n-1} g(x/|x|_2)>0.$$
By Minkowski's existence theorem \cite[p.356]{G} there is a convex origin-symmetric body $L\in
\mathbb{R}^n$ with such defined curvature function. Now multiply both sides by $|x|_2^{\alpha}$ and
apply the Fourier transform to get
\begin{eqnarray*}
-\pi (-\Delta)^{\alpha/2}\mbox{vol}_{n-1}(L|\theta^\perp)&=& - \pi(-\Delta)^{\alpha/2}\mbox{vol}_{n-1}(K|\theta^\perp)+(2\pi)^n\varepsilon  v(\theta)\\
&\ge &- \pi(-\Delta)^{\alpha/2}\mbox{vol}_{n-1}(K|\theta^\perp).
\end{eqnarray*}

On the other hand,
$$
-\pi \int_{S^{n-1}}(|x|_2^{-\alpha}
\|x\|_{K^*})^\wedge(\theta)(-\Delta)^{\alpha/2}\mbox{vol}_{n-1}(L|\theta^\perp)d\theta=$$
$$=- \pi\int_{S^{n-1}}(|x|_2^{-\alpha}
\|x\|_{K^*})^\wedge(\theta)(-\Delta)^{\alpha/2}\mbox{vol}_{n-1}(K|\theta^\perp)d\theta+$$
$$+ (2\pi)^n \varepsilon \int_{S^{n-1}}(|x|_2^{-\alpha} \|x\|_{K^*})^\wedge(\theta)   v(\theta)d\theta<$$
$$-\pi\int_{S^{n-1}}(|x|_2^{-\alpha} \|x\|_{K^*})^\wedge(\theta)(-\Delta)^{\alpha/2}\mbox{vol}_{n-1}(K|\theta^\perp)d\theta, 
$$ where the last inequality follows from the fact that $v$ is supported in the set, where
$(|x|_2^{-\alpha} \|x\|_{K^*})^\wedge<0$.

Using the argument from Theorem \ref{Thm:main2} we get that
$$\mbox{vol}_{n}(L)<\mbox{vol}_{n}(K).$$

In order to prove the remaining case when $\alpha=1$, we need two Lemmas. The following Lemma is
from \cite[Lemma 3.3]{YY}, see also \cite{KKYY}.
\begin{Lem}\label{Lem:log in R^n}
Let $K$ be an infinitely smooth origin-symmetric star body in $\mathbb{R}^n$. Then
\begin{equation} \label{logrepr in R^n}
\ln \|x\|_K =-\frac{1 }{(2\pi)^n}\int_{S^{n-1}} \ln |(x, \xi ) |  \left(\ln
\|x\|_K\right)^\wedge(\xi) d\xi + C_K,
\end{equation}
where $$C_K=\frac{1}{|S^{n-1}|}\int_{S^{n-1}} \ln\|x\|_K dx -\frac{1}{2\sqrt{\pi}}\Gamma'(1/2)+
\frac{1}{2}\frac{\Gamma'(n/2)}{\Gamma(n/2)}.$$ Moreover,
\begin{equation} \label{int_of FT of log}\int_{S^{n-1}} \left(\ln \|x\|_K\right)^\wedge(\xi) d\xi = -(2\pi)^n. \end{equation}

\end{Lem}

The following result is  from \cite[Lemma 3.7]{YY}.  It is not stated in this form there, but
follows from the proof.
\begin{Lem}\label{FT: norm to -n} Let $K$ be an origin-symmetric star body in $\mathbb{R}^n$, then
the Fourier transform of $\|x\|_K^{-n}$ is a continuous function on  $\mathbb{R}^n\setminus \{0\}$,
which equals
\begin{eqnarray}\label{FT -n} (\|x\|_K^{-n})^\wedge(\xi) &=& \int_{S^{n-1}}
\|\theta\|_K^{-n}\Big(\Gamma'(1)-\ln |(\theta,\xi)|\Big)d\theta.
\end{eqnarray}
\end{Lem}

Now we are able to prove the remaining case of the Theorem, when  $\alpha=1$ (and therefore $n=3$).
By Lemma \ref{Lem:alpha=1} there exists an infinitely smooth origin-symmetric convex body $K^\ast$
in $\mathbb{R}^3$ such that for some $\xi\in S^2$  $$ \left(|x|_2^{-1}\|x\|_{K^\ast} -
\frac{\int_{S^2} \|\theta\|_{K^\ast}d\theta }{4\pi(1+\Gamma'(1))}\ \ln |x|_2\right)^\wedge(\xi)< 0
.$$
Let $K$ be the polar body of $K^*$. By approximation we can  assume that $K$ is infinitely smooth
with strictly positive curvature,  see \cite[pp. 158-160]{S}.

 Let $$\Omega=\{\theta \in S^{2}:\left(|x|_2^{-1}\|x\|_{K^*} - \frac{\int_{S^2}
\|\theta\|_{K^*}d\theta }{4\pi(1+\Gamma'(1))}\ \ln |x|_2\right)^\wedge(\theta)< 0\}$$ and let $v\in
C^\infty (S^{2})$ be an even function, $0<v\le 1$, not identically equal to 1, and such that $ v=1$
in $S^2\setminus \Omega$. We will also use the fact that  $\Gamma'(1)>-1$ to impose an additional
condition on $v$:
$$\frac{1}{4\pi} \int_{S^2} \ln v(\theta) d\theta =
-1 - \Gamma'(1).$$  Note that the latter equality can be written in the form
\begin{equation}\label{C_V}
C_v+\Gamma'(1)=0,
\end{equation} where $C_v$ is the constant from Lemma \ref{Lem:log in R^n}.

Extend $v$ from the sphere to $\mathbb{R}^3$ as a homogeneous function of degree $1$, and denote
this extension also by $v$. By Theorem \ref{Thm:n-1-deriv} the Fourier transform of $\ln v(x) $
outside of the origin is equal to $|x|_2^{-3}g(x/|x|_2)$ for some function $g\in C^\infty(S^{2})$.
Choose an $\varepsilon>0$ small enough and define
$$ f_L(x)=f_K(x) - \varepsilon |x|_2^{-4} g(x/|x|_2)>0.$$
By Minkowski's existence theorem  there is a convex symmetric body $L\subset \mathbb{R}^3$ with
such defined curvature function. Now multiply both sides by $|x|_2$ and apply the Fourier transform
to get
\begin{eqnarray*}
&&-\pi (-\Delta)^{\alpha/2}\mbox{vol}_{2}(L|\xi^\perp)= \\
&& \hspace{2cm}= - \pi(-\Delta)^{\alpha/2}\mbox{vol}_{2}(K|\xi^\perp)-(2\pi)^3\varepsilon (|x|_2^{-3} g(x/|x|_2))^\wedge(\xi)\\
&&\hspace{2cm}\ge - \pi(-\Delta)^{\alpha/2}\mbox{vol}_{2}(K|\xi^\perp),
\end{eqnarray*}
where the last inequality comes from the following calculations, based on Lemmas \ref{FT: norm to
-n} and \ref{Lem:log in R^n}.
$$\left(|x|_2^{-3} g(x/|x|_2)\right)^\wedge(\xi) = \int_{S^{2}} g(\theta) (\Gamma'(1) - \ln |(\theta,\xi)|)
d\theta$$
$$= \int_{S^{2}} \left(\ln v(x)\right)^\wedge(\theta) (\Gamma'(1) - \ln |(\theta,\xi)|)
d\theta$$
$$= -(2\pi)^3 \Gamma'(1) - \int_{S^{2}} \left(\ln v(x) \right)^\wedge(\theta)  \ln |(\theta,\xi)|
d\theta$$
$$= -(2\pi)^3 \Gamma'(1) +(2\pi)^3 ( \ln v(\xi) - C_v) =  (2\pi)^3  \ln v(\xi)\le 0.$$

On the other hand
$$
  V_1(L,K)= \frac13 \int_{S^{2}}\|\theta\|_{K^*}f_L(\theta) d\theta =
\frac13\int_{S^{2}}\|\theta\|_{K^*}(f_K(\theta) - \varepsilon  g(\theta))d\theta$$
$$= \mbox{vol}_{3}(K) - \frac{\varepsilon}{3} \int_{S^{2}}\|\theta\|_{K^*}  g(\theta)d\theta.
$$

If we can show that $\int_{S^{2}}\|\theta\|_{K^*}  g(\theta)d\theta> 0$, the statement will follow
from Minkowski's first inequality. Let  $\|\theta\|_M= \exp{\|\theta\|_{K^*}}$ for all $\theta\in
S^{2}$. 
By Lemma \ref{Lem:log in
R^n} we have
$$\|\theta\|_{K^*} =  \ln \|\theta\|_M= -\frac{1}{(2\pi)^3}\int_{S^{2}} \ln |(\theta, \xi)| (\ln \|x\|_M)^\wedge(\xi) d\xi+ C_M,$$
where
 \begin{equation}\label{C_M}
 C_M= \frac{1}{4\pi}\int_{S^2} \ln \|\theta\|_M d\theta+1=\frac{1}{4\pi}\int_{S^2} \|\theta\|_{K^*} d\theta+1.
 \end{equation}
Analogously,
$$ \ln |\theta|_2 = -\frac{1}{(2\pi)^3}\int_{S^{2}} \ln |(\theta, \xi)| (\ln |x|_2)^\wedge(\xi) d\xi+ 1.$$
Let us denote $$\lambda=\frac{C_M+\Gamma'(1)}{1+\Gamma'(1)}.$$ Then
$$\int_{S^{2}}\|\theta\|_{K^*}  g(\theta)d\theta= \int_{S^{2}}\left( \ln \|\theta\|_{M} - \lambda \ln |\theta |_2 \right)g(\theta)  d\theta=
$$
$$=-\frac{1}{(2\pi)^3}\int_{S^{2}}\left( \int_{S^{2}}\ln |(\theta, \xi)|  \left(\ln \|x\|_M - \lambda \ln |x|_2\right)^\wedge(\xi) d\xi\right)
g(\theta)d\theta$$ $$  +  \left(C_M - \lambda \right) \int_{S^{2}} g(\theta)d\theta.$$ Reversing
the order of integration in the first integral and  then adding and subtracting an appropriate
quantity, we get

$$=\frac{1}{(2\pi)^3}\int_{S^{2}}\left( \int_{S^{2}}\left(\Gamma'(1) - \ln |(\theta, \xi)| \right) g(\theta)d\theta \right) \left(\ln \|x\|_M - \lambda \ln |x|_2\right)^\wedge(\xi) d\xi
$$
$$ - \frac{\Gamma'(1)}{(2\pi)^3}\int_{S^{2}}\left( \int_{S^{2}}  g(\theta)d\theta \right) \left(\ln \|x\|_M - \lambda \ln |x|_2\right)^\wedge(\xi) d\xi $$
$$ +  \left(C_M - \lambda \right) \int_{S^{2}} g(\theta)d\theta.$$
Formulas (\ref{FT -n}) and (\ref{int_of FT of log})  applied to the first and second integrals
correspondingly give
$$=\frac{1}{(2\pi)^3}\int_{S^{2}}\left(|x|_2^{-3} g(x/|x|_2)\right)^\wedge(\xi) \left(\ln \|x\|_M - \lambda \ln |x|_2\right)^\wedge(\xi) d\xi
$$
$$ +{\Gamma'(1)} (1- \lambda)   \int_{S^{2}}  g(\theta)d\theta  +  \left(C_M - \lambda \right) \int_{S^{2}} g(\theta)d\theta.$$
Using that $$\left(|x|_2^{-3} g(x/|x|_2)\right)^\wedge(\xi)=(2\pi)^3  \ln v(\xi) $$and
$${\Gamma'(1)} (1- \lambda) + \left(C_M - \lambda \right)=0,$$
 we   get
\begin{equation}\label{eqn:final}
\int_{S^{2}}\|\theta\|_{K^*}  g(\theta)d\theta=\int_{S^{2}}\ln v(\xi) \left(\ln \|x\|_M - \lambda
\ln |x|_2\right)^\wedge(\xi) d\xi. \end{equation}
 Recall that
$$ \ln \|x\|_M - \lambda \ln |x|_2 $$
$$= |x|_2^{-1}\|x\|_{K^*}+\ln |x|_2 - \frac{\frac{1}{4\pi}\int_{S^2} \|\theta\|_{K^*}
d\theta+1+\Gamma'(1)}{1+\Gamma'(1)} \ln|x|_2$$
$$=|x|_2^{-1}\|x\|_{K^\ast} -
\frac{\int_{S^2} \|\theta\|_{K^\ast}d\theta }{4\pi(1+\Gamma'(1))}\ \ln |x|_2.$$ Therefore
(\ref{eqn:final}) implies
$$\int_{S^{2}}\|\theta\|_{K^*}  g(\theta)d\theta>0,$$
 since $\ln v$ is negative, where $$\left(|x|_2^{-1}\|x\|_{K^\ast} -
\frac{\int_{S^2} \|\theta\|_{K^\ast}d\theta }{4\pi(1+\Gamma'(1))}\ \ln |x|_2\right)^\wedge$$ is
negative, and zero everywhere else.

 \qed

\noindent{\bf Remark.} The aim of Theorem \ref{Thm:main3} is to show that condition
(\ref{eqn:Laplacians}) is inconclusive. As we have seen, there are  bodies for which
(\ref{eqn:Laplacians}) holds, but $\mathrm{vol}_n(L)< \mathrm{vol}_n(K).$ Let us remark that one
can easily find two bodies for which (\ref{eqn:Laplacians}) holds, but $\mathrm{vol}_n(L)>
\mathrm{vol}_n(K).$ It is enough to take two Euclidean balls. This is obvious for $\alpha> 1$, but
probably some explanations are needed for the case $\alpha=1$ (and  $n=3$).

Let $B$ be a Euclidean ball in $\mathbb{R}^3$ with curvature function $f_B(\theta) =C$, $\forall
\theta\in S^{2}$. Then by (\ref{eqn:Def-Laplace}) and (\ref{FT -n}) we have
$$(-\Delta)^{1/2}\mathrm{vol}_{2}(B|\xi^\perp)= -\frac{1}{\pi}( |x|_2
f_B(x))^\wedge(\theta)= -\frac{1}{\pi}( |x|_2  C |x|_2^{-4})^\wedge(\theta)=$$
$$=-\frac{C}{\pi}\int_{S^{2}}  \Big(\Gamma'(1)-\ln |(\theta,\xi)|\Big)d\xi.$$
The latter integral is computable and after routine calculations one gets
$$(-\Delta)^{1/2}\mathrm{vol}_{2}(B|\xi^\perp) = - 4 C(\Gamma'(1)+1)<0.$$
Therefore if we take two Euclidean balls $B_r$ and $B_R$ with radii $r<R$, then $f_{B_r}<f_{B_R}$
and therefore
$$ (-\Delta)^{1/2}\mathrm{vol}_{2}(B_R|\xi^\perp) < (-\Delta)^{1/2}\mathrm{vol}_{2}(B_r|\xi^\perp),$$
but $$\mathrm{vol}_3(B_r)< \mathrm{vol}_3(B_R).$$


\begin{thebibliography}{WWW}


\bibitem[B]{B} {\sc K. Ball}, {\em Shadows of convex bodies}, Trans. Amer. Math. Soc. {\bf 327} (1991) ,
891--901.

\bibitem[G]{G} {\sc R.~J.~Gardner}, {\em Geometric Tomography}, Cambridge
University Press, Cambridge 1995.

\bibitem[GKS]{GKS}
{\sc R.~J.~Gardner, A.~Koldobsky, T.~Schlumprecht}, {\em An analytic solution to the Busemann-Petty
problem on sections of convex bodies}, Annals of Math. {\bf 149} (1999), 691--703.

\bibitem[GZ]{GZ} {\sc P.~Goodey and G.~Zhang}, {\em Inequalities between projection functions of
convex bodies}, Amer. J. Math. {\bf 120} (1998), 345--367.

\bibitem[K1]{K-Israel99}
{\sc A.~Koldobsky},  {\em A generalization of the Busemann-Petty problem on sections of convex
bodies}, Israel J. Math. {\bf 110} (1999), 75--91.

\bibitem[K2]{K-AdvAppl}
{\sc A.~Koldobsky}, {\em Comparison of volumes by means of the areas of central sections},  Adv. in
Appl. Math.  {\bf 33}  (2004),  no. 4, 728--732.

\bibitem[K3]{K-book}
{\sc A.~Koldobsky,} {\em Fourier Analysis in Convex Geometry}, Mathematical Surveys and Monographs,
American Mathematical Society, Providence RI, 2005.


\bibitem[KKYY]{KKYY} {\sc N.~J.~Kalton, A.~Koldobsky, V.~Yaskin and
M.~Yaskina,} {\em The geometry of $L_0$}, Canadian J. Math., to appear.
http://arxiv.org/abs/math/0412371


\bibitem[KRZ1]{KRZ1}{\sc A.~Koldobsky, D.~Ryabogin and  A.~Zvavitch,}
{\em Projections of convex bodies and the Fourier transform},  { Israel J. Math.} {\bf  139}
(2004), 361--380.

\bibitem[KRZ2]{KRZ2}{\sc A.~Koldobsky, D.~Ryabogin and  A.~Zvavitch,}
{\em Unified Fourier analytic approach to the volume of projections and sections of convex bodies},
Fourier Analysis and Convexity, (Editors: L. Brandolini, L. Cozani, A. Iosevich and G. Travaglini),
Birkhauser 2004, 119-131.

\bibitem[KYY]{KYY} {\sc A.~Koldobsky, V.~Yaskin and M.~Yaskina}, {\em Modified Busemann-Petty problem on sections on
convex bodies}, Israel J. Math., {\bf 154} (2006), 191--208.

\bibitem[P]{Pe}{\sc C.~M.~Petty}, {\em Projection bodies}, {
Proc.~Coll.~Convexity (Copenhagen 1965), Kobenhavns Univ.~Mat.~Inst.,} 234--241.


\bibitem[RZ]{RZ}
{\sc D.~Ryabogin,  A.~Zvavitch}, {\em The Fourier transform and Firey projections of convex
bodies}, Indiana Univ. Math. J.  {\bf 53}  (2004),  no. 3, 667--682.


\bibitem[S1]{S1} {\sc R.~Schneider}, {\em Zur einem Problem von Shephard \"uber die
Projektionen konvexer K\"orper}, { Math.~Z.}~{\bf 101} (1967), 71--82.

\bibitem[S2]{S}{\sc R.~Schneider},
{\em Convex Bodies: the Brunn-Minkowski Theory}, Cambridge University Press, Cambridge, 1993.


\bibitem[YY]{YY} {\sc  V.~Yaskin and M.~Yaskina,} {\em Centroid bodies and
comparison of volumes}, Indiana University Math. J.,  {\bf 55} No. 3 (2006), 1175–-1194.



\end{thebibliography}
\end{document}